\documentclass[12pt]{amsart}
\usepackage{amsmath,amssymb,amsbsy,amsfonts,amsthm,latexsym,amsopn,amstext,
                                                       amsxtra,euscript,amscd}
\begin{document}

\title[Formulas for Positive, Negative and Zero Values \ldots]{Formulas for Positive, Negative and Zero Values of the M\"obius Function}

\author[R. M.~Abrarov]{R. M.~Abrarov$^1$}
%\address{University of Toronto, Canada}
%\rabrarov@physics.utoronto.ca
\email{$^1$rabrarov@gmail.com}
\author[S. M.~Abrarov]{S. M.~Abrarov$^2$}
%\address{York University, Toronto, Canada}
%\email{abrarov@yorku.ca}
\email{$^2$absanj@gmail.com}

\date{May 4, 2009}

\begin{abstract}

We obtained the formulas for the quantities of positive, negative and zero values of the M\"obius function for any real $x$ in terms of the M\"obius function values for square root of $x$ -- similar to the identities we found earlier for the Mertens function [1]. Using the Dirac delta function approach [3, 2] we propose the equations showing how the Mertens and related functions can generate the output values.

\end{abstract}

\maketitle

\section{\textbf{Distributions of the M\"obius function values}}

The Mertens function $M(x)=\sum _{k=1}^{x}\mu _{k}$ can be expressed through identity \cite{RAbrarov}
\begin{equation} \label{Eq_1}
M\left(x\right)=2M\left(\sqrt{x} \right)-\sum _{i,j=1}^{\sqrt{x} }\mu _{i}  \mu _{j} \left[\frac{x}{i\cdot j} \right],     
\end{equation}
where $\mu _{k} $ is the M\"obius function. It is possible to find similar identities for distributions of positive, negative and zero values of the M\"obius function. 

Let us denote $N_{+} (x),\, N_{-} (x),\, N_{0} (x)$ as the quantity of natural numbers correspondingly with $\mu =+1$, $\mu =-1$ and $\, \mu =0$ in the interval not exceeding $x\ge 0$. For the Mertens function we have $M(x)=N_{+} (x)-N_{-} (x)$. We will also apply a notation $Q(x)$ \cite{Hardy, Weisstein1}, commonly used for the quantity of squarefree numbers $N_{+} (x)+N_{-} (x)$.

There are obvious relations between quantities
\begin{equation} \label{Eq_2} 
\begin{array}{l} {\left[x\right]=N_{+} (x)+\, N_{-} (x)+N_{0} (x)\, ,} \\ {2N_{+} (x)=\left[x\right]-N_{0} (x)+M(x)\, ,} \\ {2N_{-} (x)=\left[x\right]-N_{0} (x)-M(x)\, .} \end{array} 
\end{equation} 

The frequencies of appearing numbers with $\mu =+1$, $\mu =-1$ and $\mu =0$ are $\frac{N_{+} (x)}{\left[x\right]} $, $\frac{N_{-} (x)}{\left[x\right]} $ and $\frac{N_{0} (x)}{\left[x\right]} $, respectively. It should be noted that the frequencies have some properties, similar to properties of the probabilities such as
\begin{equation} \label{Eq_3} 
\frac{N_{+} (x)}{\left[x\right]} +\frac{N_{-} (x)}{\left[x\right]} +\, \frac{N_{0} (x)}{\left[x\right]} =1.      
\end{equation} 
However we are considering here the deterministic process. Therefore these frequencies must not be confused with probabilities. 

Let us find identities for $N_{+} (x),\, N_{-} (x),\, N_{0} (x)$ through known distribution of the M\"obius function for $\sqrt{x} $, similar to \eqref{Eq_1}. It is straightforward to get expression for $N_{0} (x)$ by using the \textit{inclusion-exclusion principle} \cite{Weisstein2, Havil}. In this case we count the quantity of numbers divisible by $2^{2} $ plus numbers divisible by $3^{2} $, except those divisible by $2^{2} $ and $3^{2} $, and so on for squares of all prime numbers:
\begin{equation} \label{Eq_4} 
\begin{aligned}
N_{0} (x)=
& \left[\frac{x}{2^{2} } \right]+\left[\frac{x}{3^{2} } \right]-\left[\frac{x}{2^{2} \cdot 3^{2} } \right] \\
&+ \left[\frac{x}{5^{2} } \right]-\left[\frac{x}{2^{2} \cdot 5^{2} } \right]-\left[\frac{x}{3^{2} \cdot 5^{2} } \right]+\left[\frac{x}{2^{2} \cdot 3^{2} \cdot 5^{2} } \right]+\ldots \,.
\end{aligned}   
\end{equation} 

If we rearrange the terms in order of the appearance of nonzero elements, we will get
\begin{equation} \label{Eq_5} 
N_{0} \left(x\right)=-\sum _{i=2}^{\sqrt{x} }\mu _{i}  \left[\frac{x}{i^{2} } \right] .       
\end{equation}
Hence from \eqref{Eq_5}, \eqref{Eq_1} and \eqref{Eq_2} immediately follow
\begin{equation} \label{Eq_6} 
\begin{aligned}
& Q\left(x\right)=\sum _{i=1}^{\sqrt{x} }\mu _{i}  \left[\frac{x}{i^{2} } \right]\, , \\
& N_{+} \left(x\right)=M\left(\sqrt{x} \right)-\frac{1}{2} \sum _{i,j=1}^{\sqrt{x} }\mu _{i}  \mu _{j} \left[\frac{x}{i\cdot j} \right]+\frac{1}{2} \sum _{i=1}^{\sqrt{x} }\mu _{i}  \left[\frac{x}{i^{2} } \right]\, , \\ 
& N_{-} \left(x\right)=-M\left(\sqrt{x} \right)+\frac{1}{2} \sum _{i,j=1}^{\sqrt{x} }\mu _{i}  \mu _{j} \left[\frac{x}{i\cdot j} \right]+\frac{1}{2} \sum _{i=1}^{\sqrt{x} }\mu _{i}  \left[\frac{x}{i^{2} } \right]\, .
\end{aligned}
\end{equation}

\section{\textbf{Dirac and Kronecker delta function transformation}}

Formula \eqref{Eq_1} involves the floor function and explicit dependence on the Mertens function in the right hand side. In order to avoid them, we can apply an alternative for the floor function - the Dirac delta function $\delta ^{D} $ \cite{Antosik} (see also our \cite{SAbrarov}) or the discrete version of the Dirac delta function $\delta ^{K} (x)\equiv \delta _{\left\lfloor x\right\rfloor \left\lceil x\right\rceil } $ -- the Kronecker delta with floor and ceiling functions of $x$ as arguments (which means that the value of $\delta ^{K} (x)$ is 1 only for integer \textit{x} and 0 otherwise)
\begin{equation} \label{Eq_7} 
\begin{aligned}
\left[\frac{x}{i\cdot j} \right] 
&=\int _{i\cdot j-\varepsilon }^{x}\delta ^{D} \left(\frac{y}{i\cdot j} -\left[\frac{y}{i\cdot j} \right]\right)dy \\
&=\int _{1-\varepsilon }^{x}\delta ^{D} \left(\frac{y}{i\cdot j} -\left[\frac{y}{i\cdot j} \right]\right)dy ,\, \, \, \, \, \, {\rm here\; and\; further}\, x\ge 1\, , 
\end{aligned}
\end{equation} 
\begin{equation} \label{Eq_8} 
\left[\frac{x}{i\cdot j} \right]=\sum _{k=i\cdot j}^{x}\delta ^{K} \left(\frac{k}{i\cdot j} \right) =\sum _{k=1}^{x}\delta ^{K} \left(\frac{k}{i\cdot j} \right)  .        
\end{equation} 
Now substituting \eqref{Eq_7} and \eqref{Eq_8} into \eqref{Eq_1}, we derive

\begin{equation} \label{Eq_9} 
\begin{aligned}
M(x)=2 &-\int _{1^{2} -\varepsilon }^{x\ge 1^{2} }\mu _{1} \left(\mu _{1} \delta ^{D} \left(\frac{y}{1\cdot 1} -\left[\frac{y}{1\cdot 1} \right]\right)\right) dy \\ 
&-\int _{2^{2} -\varepsilon }^{x\ge 2^{2} }\mu _{2} \left(2\mu _{1} \delta ^{D} \left(\frac{y}{2\cdot 1} -\left[\frac{y}{2\cdot 1} \right]\right) \right. \\
& +\mu _{2} \delta ^{D} \left. \left(\frac{y}{2\cdot 2} -\left[\frac{y}{2\cdot 2} \right]\right)\right) dy \\
&-\int _{3^{2} -\varepsilon }^{x\ge 3^{2} }\mu _{3} \left(2\mu _{1} \delta ^{D} \left(\frac{y}{3\cdot 1} -\left[\frac{y}{3\cdot 1} \right]\right)\right. \\
&+2\mu _{2} \delta ^{D} \left. \left(\frac{y}{3\cdot 2} -\left[\frac{y}{3\cdot 2} \right]\right)+\mu _{3} \delta ^{D} \left(\frac{y}{3\cdot 3} -\left[\frac{y}{3\cdot 3} \right]\right)\right) dy \\ 
&-\int _{5^{2} -\varepsilon }^{x\ge 5^{2} }\mu _{5} \left(2\mu _{1} \delta ^{D} \left(\frac{y}{5\cdot 1} -\left[\frac{y}{5\cdot 1} \right]\right) \right. \\
&+2\mu _{2} \delta ^{D} \left(\frac{y}{5\cdot 2} -\left[\frac{y}{5\cdot 2} \right]\right)+2\mu _{3} \delta ^{D} \left(\frac{y}{5\cdot 3} -\left[\frac{y}{5\cdot 3} \right]\right)\\
&+\mu _{5} \delta ^{D} \left. \left(\frac{y}{5\cdot 5} -\left[\frac{y}{5\cdot 5} \right]\right)\right) dy-\ldots\\
=2 &-\sum _{i=1}^{\sqrt{x} }\int _{i^{2} -\varepsilon }^{x}\mu _{i} \left(-\mu _{i} \delta ^{D} \left(\frac{y}{i\cdot i} -\left[\frac{y}{i\cdot i} \right]\right)\right. \\
&+\sum _{j=1}^{i}2\mu _{j}  \delta ^{D} \left.\left(\frac{y}{i\cdot j} -\left[\frac{y}{i\cdot j} \right]\right)\right) dy ,
\end{aligned} 
\end{equation}
\begin{equation} \label{Eq_10} 
\begin{aligned}
M(x)& =2-\sum _{k=1}^{x}\sum _{i,j=1}^{\sqrt{k} }\mu _{i}   \mu _{j} \delta ^{K} \left(\frac{k}{i\cdot j} \right) \\
& =2-\sum _{i=1}^{\sqrt{x} }\sum _{k=i^{2} }^{x}\mu _{i}   \left(-\mu _{i} \delta ^{K} \left(\frac{k}{i\cdot i} \right)+\sum _{j=1}^{i}2\mu _{j}  \delta ^{K} \left(\frac{k}{i\cdot j} \right)\right).
\end{aligned}    
\end{equation}
From \eqref{Eq_10} we have 

\begin{equation} \label{Eq_11} 
\mu _{x} =-\sum _{i,j=1}^{\sqrt{x} }\mu _{i}  \mu _{j} \delta ^{K} \left(\frac{x}{i\cdot j} \right)\qquad \qquad x\ge 2.      
\end{equation}

By following the same procedure as it was in the previous section, we obtain

\begin{equation} \label{Eq_12} 
\begin{aligned}
N_{0} \left(x\right) = & -\sum _{i=2}^{\sqrt{x} }\mu _{i} \int _{i^{2} -\varepsilon }^{x}\delta ^{D} \left(\frac{y}{i\cdot i} -\left[\frac{y}{i\cdot i} \right]\right)  dy \\
= & -\sum _{i=2}^{\sqrt{x} }\sum _{k=i^{2} }^{x}\mu _{i} \delta ^{K} \left(\frac{k}{i\cdot i} \right) ,
\end{aligned}
\end{equation}
\begin{equation} \label{Eq_13} 
\begin{aligned} 
Q\left(x\right) = & \sum _{i=1}^{\sqrt{x} }\mu _{i} \int _{i^{2} -\varepsilon }^{x}\delta ^{D} \left(\frac{y}{i\cdot i} -\left[\frac{y}{i\cdot i} \right]\right)  dy = \sum _{i=1}^{\sqrt{x} }\sum _{k=i^{2} }^{x}\mu _{i} \delta ^{K} \left(\frac{k}{i\cdot i} \right) ,
\end{aligned}
\end{equation}
\begin{equation} \label{Eq_14} 
\begin{aligned}
N_{+} \left(x\right) = & 1-\frac{1}{2} \sum _{i=1}^{\sqrt{x} }\int _{i^{2} -\varepsilon }^{x}\mu _{i} \left(-\mu _{i} \delta ^{D} \left(\frac{y}{i\cdot i} -\left[\frac{y}{i\cdot i} \right]\right)\right. \\
& \,\, +\sum _{j=1}^{i}2\mu _{j}  \delta ^{D} \left. \left(\frac{y}{i\cdot j} -\left[\frac{y}{i\cdot j} \right]\right)\right)  dy \\ 
& \,\, +\frac{1}{2} \sum _{i=1}^{\sqrt{x} }\mu _{i} \int _{i^{2} -\varepsilon }^{x}\delta ^{D} \left(\frac{y}{i\cdot i} -\left[\frac{y}{i\cdot i} \right]\right)  dy \\ 
 = & 1-\frac{1}{2} \sum _{i=1}^{\sqrt{x} }\sum _{k=i^{2} }^{x}\mu _{i}   \left(-\mu _{i} \delta ^{K} \left(\frac{k}{i\cdot i} \right) \right. \\
& \,\, +\sum _{j=1}^{i}2\mu _{j}  \delta ^{K} \left. \left(\frac{k}{i\cdot j} \right)\right) +\frac{1}{2} \sum _{i=1}^{\sqrt{x} }\sum _{k=i^{2} }^{x}\mu _{i} \delta ^{K} \left(\frac{k}{i\cdot i} \right) ,
\end{aligned}
\end{equation}
\begin{equation} \label{Eq_15} 
\begin{aligned}
N_{-} \left(x\right) = & -1+\frac{1}{2} \sum _{i=1}^{\sqrt{x} }\int _{i^{2} -\varepsilon }^{x}\mu _{i} \left(-\mu _{i} \delta ^{D} \left(\frac{y}{i\cdot i} -\left[\frac{y}{i\cdot i} \right]\right)\right. \\
& +\sum _{j=1}^{i}2\mu _{j}  \delta ^{D} \left. \left(\frac{y}{i\cdot j} -\left[\frac{y}{i\cdot j} \right]\right)\right)  dy \\
& +\frac{1}{2} \sum _{i=1}^{\sqrt{x} }\mu _{i} \int _{i^{2} -\varepsilon }^{x}\delta ^{D} \left(\frac{y}{i\cdot i} -\left[\frac{y}{i\cdot i} \right]\right)  dy \\ 
= & -1+\frac{1}{2} \sum _{i=1}^{\sqrt{x} }\sum _{k=i^{2} }^{x}\mu _{i}   \left(-\mu _{i} \delta ^{K} \left(\frac{k}{i\cdot i} \right)\right. \\
& +\sum _{j=1}^{i}2\mu _{j}  \delta ^{K} \left. \left(\frac{k}{i\cdot j} \right)\right) +\frac{1}{2} \sum _{i=1}^{\sqrt{x} }\sum _{k=i^{2} }^{x}\mu _{i} \delta ^{K} \left(\frac{k}{i\cdot i} \right) . 
\end{aligned}
\end{equation}

These identities are showing a clear mechanism how functions above generate the output values and can be helpful in understanding of their distributions and other properties.
\\
\\

\end{document}